\documentstyle{amsppt}
\magnification=1200
\NoRunningHeads

\define\0{\bold0}
\define\1{\bold1}
\define\RR{\Bbb R}
\define\Hom{\operatorname{Hom}}
\define\Low{\operatorname{Low}}
\define\sT{\Cal T}
\define\sB{\Cal B}
\let\emptyset\varnothing
\let\phi\varphi
\let\epsilon\varepsilon
\let\wt\widetilde
\let\wh\widehat

\topmatter
\title{Tensor Products of Idempotent Semimodules.\\
An Algebraic Approach}\endtitle

\author{G.~L.~Litvinov, V.~P.~Maslov, and G.~B.~Shpiz}\endauthor

\thanks
Appears in: {\it Math. Notes\/} {\bf 64}:4 (1999) 479--489, as the English
translation of the paper published in {\it Mat.  Zametki\/} {\bf 64}:4
(1999) 572--585.
\endthanks


%
\abstract\nofrills
{\smc Abstract}. 
We study idempotent analogs of topological tensor
products in the sense of A. Grothendieck. The basic concepts and
results are simulated on the algebraic level. This is one of a
series of papers on idempotent functional analysis.
\newline\newline\noindent
{\smc Key words}: \, 
idempotent functional analysis, idempotent semiring, idempotent
semimodule, tensor product, polylinear mapping, nuclear operator.
\endabstract
\endtopmatter

\document
\rightline{\it Dedicated to S.G. Krein on the occasion of his 80th birthday}
\head
Introduction
\endhead

We construct and study tensor products in some natural
categories of idempotent semimodules. In idempotent analysis,
these tensor products play a role similar to that of the
topological tensor products constructed by
A. Grothendieck~\cite{1} in functional analysis. However, we
point out that the idempotent versions of the basic notions and
results are fairly different from their conventional analogs.

The basic concepts and results (including those of
``topological'' nature) are simulated on the algebraic
level; the point is that the operation of idempotent
addition can be defined for infinitely many
summands. The present paper is one of a series of papers
on idempotent functional analysis, which is an
abstract version of idempotent analysis in the sense
of~\cite{2--9}.
In the subsequent publications, we intend to study
the links between idempotent tensor
products, idempotent linear integral and
nuclear operators, and traces of operators.

\head
\S1. Basic Concepts
\endhead

\subhead
1.1
\endsubhead
We recall that an {\it idempotent semigroup\/}
is an additive semigroup $S$ with
commutative addition $\oplus$ such that
$x\oplus x = x$ for all elements $x\in S$.
If $S$ has a neutral element, then this element
is denoted by $\0$. Any idempotent semigroup is
a partially
ordered set with respect to the {\it standard
order\/} defined as follows: $x\preccurlyeq y$ if and only if
$x\oplus y = y$.
This is obviously well defined, and
$x\oplus y = \sup\{x,y\}$.
Thus each idempotent semigroup is an
upper semilattice, and moreover, the notions of an
idempotent semigroup and an upper
semilattice essentially coincide~\cite{10}.

The definitions given below are partly borrowed from~\cite{11,~12}.
An idempotent semigroup $S$
is said to be $a$-{\it complete\/} (or {\it algebraically
complete\/}) if it is complete as an ordered set,
i.e., each subset $X\subset S$ has the least upper
bound $\sup X$,
denoted also by $\oplus X$, and the greatest
lower bound $\inf X$, denoted also by $\wedge X$.
This semigroup is said to be $b$-{\it complete\/} (or
{\it boundedly complete\/}) if each subset $X\subset S$
(possibly, empty) that is bounded above has the least
upper bound $\oplus X$ (in this case each nonempty
subset $Y\subset S$ has the greatest lower bound
$\wedge Y$, and $S$ is a lattice). Note that each
$a$- or $b$-complete idempotent semigroup
has a zero $\0$, which coincides with
$\oplus\emptyset$, where $\emptyset$ is the empty set.
Needless to say, $a$-completeness implies $b$-completeness.
The cut completion procedure~\cite{10}
provides an embedding $S\to\wh S$
of an arbitrary idempotent semigroup $S$ in an
$a$-complete idempotent semigroup $\wh S$
(called the {\it normal completion\/} of $S$); furthermore,
$\widehat{\widehat S\,}\! = S$. In a similar way, one defines the
$b$-completion $S\to\wh S_b$: if $S\ni\infty = \sup S$, then
$\wh S_b = S$; otherwise, $\wh S = \wh S_b\cup\{\infty\}$.
An arbitrary $b$-complete idempotent semigroup $S$ also can
differ from $\wh S$ only by the element $\infty = \sup S$.

Let $S$ and $T$ be $b$-complete idempotent semigroups. A
homomorphism $f\:S\to T$ will be called a $b$-{\it
homomorphism\/}
if $f(\oplus X) = \oplus f(X)$ for each bounded
subset $X$ of $S$. If a $b$-homomorphism $f$
extends to be a homomorphism $\wh S\to\wh T$ of the normal
completions, and moreover, if $f(\oplus X) = \oplus f(X)$ for all
$X\subset S$, then $f$ will be called an $a$-{\it homomorphism}.
Each $a$-homomorphism is a $b$-homomorphism. If the
semigroup $S$ is $a$-complete, then each
$b$-homomorphism is an $a$-homomorphism. If $S$ and $T$
are topological idempotent semigroups, then a homomorphism $f\:T\to S$
that takes zero to zero is an $a$-homomorphism if and only if
$f$ is lower semicontinuous~\cite{11,~12}.

\subhead
1.2
\endsubhead
An {\it idempotent semiring\/} (for brevity, we sometimes
simply say ``semiring'') is an idempotent
semigroup $(K,\oplus)$ equipped with an associative
multiplication $\odot$ such that both distributivity
conditions hold. If the multiplication is commutative, then
the idempotent semiring
is said to be {\it commutative}. An element $\1\in K$
is called the {\it unit\/} of the semirings $K$ if
$x\odot\1=\1\odot x = x$ for all $x\in K$. An element
$\0\in K$, $\0\ne\1$, is called the {\it zero\/}
of the semirings $K$ if $x\oplus\0 = x$ and
$x\odot\0=\0\odot x=\0$ for all $x\in K$. In this paper we
consider only idempotent semirings with unit. The presence of
zero is usually also assumed (unless specified otherwise).
An {\it idempotent semifield\/} (or simply
{\it semifield\/}) is a commutative idempotent semiring in
which every nonzero element is invertible with respect to
multiplication.
An idempotent semiring $K$ is said to be $a$-{\it complete\/}
(respectively, $b$-{\it complete\/}) if $K$ is an $a$-complete
(respectively, $b$-complete) idempotent semigroup and if
the generalized distributivity laws $k\odot(\oplus X)=\oplus(k\odot X)$
and $(\oplus X)\odot k=\oplus(X\odot k)$ hold
for each subset
(respectively, each bounded subset) $X\subset K$ and each $k\in K$.
It follows from the generalized distributivity that each $a$-complete or
$b$-complete idempotent semiring has a zero element, which
coincides with
$\oplus\emptyset$. The notion of an $a$-complete idempotent semiring
coincides with that of a complete dioid in the sense of~\cite{13,~14}.

The set $\RR(\max,+)$ of real numbers equipped with the
idempotent addition $\oplus = \max$ and the multiplication
$\odot= +$ is an example of an idempotent semiring; in this case
$\1= 0$.
By supplementing this idempotent semiring with the
element $\0 = -\infty$, we obtain the
$b$-complete semiring $\RR_{\max} = \RR\cup\{-\infty\}$ with the
same operations and with zero. By supplementing $\RR_{\max}$
with the element $+\infty$ and by setting $\0\odot(+\infty) = \0$,
$x\odot(+\infty) = +\infty$ for $x\ne\0$,
$x\oplus(+\infty) = +\infty$ for all $x$, we obtain the
$a$-complete idempotent semiring
$\wh{\RR}_{\max} = \RR_{\max}\cup\{+\infty\}$. The standard
order on $\RR(\max,+)$, $\RR_{\max}$, and $\wh{\RR}_{\max}$
coincides with the usual order. The idempotent semirings $\RR(\max,+)$,
$\RR_{\max}$ are semifields. On the other hand, an
$a$-complete idempotent semiring other than $\{\0,\1\}$ cannot
be a semifield. There is an important class of examples
related to (topological) vector lattices (e.g., see~\cite{15,
Chap.~V} and~\cite{10}).
By defining the sum $x\oplus y$
as $\sup\{x,y\}$ and the multiplication $\odot$ as addition
of vectors, one can interpret vector lattices as
idempotent semifields. By supplementing a complete vector
lattice (in the sense of~\cite{10,~15})
with the element $\0$,
we obtain a $b$-complete semifield. Next, we can add the
``infinite'' element, thus obtaining an $a$-complete
idempotent semiring (coinciding as the
ordered set with the normal completion of the original
lattice).

\subhead
1.3
\endsubhead
Let $V$ be an idempotent semigroup and $K$ an
idempotent semiring. Suppose that a multiplication
$k,x\mapsto k\odot x$ of elements of $V$ by elements of $K$ is
given,
and moreover, this multiplication is associative, distributes
over the addition in $V$, and satisfies $\1\odot x= x$ for all
$x\in V$. In this case the semigroup $V$ is called an
{\it idempotent semimodule\/} (or simply
{\it semimodule\/}) over $K$. An element $\0_V\in V$
is called the {\it zero of the semimodule\/} $V$ if
$k\odot\0_V = \0_V$, $\0_K\odot x = \0_V$, and $\0_V\oplus x = x$
for all $k\in K$ and $x\in V$. Let $V$ be a semimodule over a
$b$-complete idempotent semiring $K$. This semimodule
is said to be $b$-{\it complete} if it is $b$-complete as an
idempotent semigroup and if
the generalized distributivity laws
$(\oplus Q)\odot x = \oplus(Q\odot x)$ and
$k\odot(\oplus X) = \oplus(k\odot X)$, $k\in K$, $x\in X$,
hold for any bounded subsets $Q\subset K$ and $X\subset V$.
This semimodule is said to be $a$-{\it complete} if
it is $b$-complete and contains the element $\infty = \sup V$.

Let $V$ and $W$ be idempotent semimodules over an
idempotent semiring $K$. A mapping $p\:V\to W$
is said to be {\it linear\/} (over $K$) if
$p(x\oplus y) = p(x)\oplus p(y)$ and $p(k\odot x) = k\odot p(x)$
for all $x,y\in V$ and $k\in K$. Let the semimodules $V$ and $W$
be $b$-complete. A linear mapping $p\:V\to W$ is said to be
$b$-{\it linear\/} if it is a $b$-homomorphism of idempotent
semigroups and $a$-{\it linear\/} if it
extends to be an $a$-homomorphism of the normal
completions $\wh V$ and $\wh W$. Every $a$-linear
mapping is $b$-linear. If a mapping
$V\to W$ is $b$-linear and the semimodule $V$ is
$a$-complete, then this mapping is $a$-linear.
The $a$-linearity is similar to continuity
(semicontinuity) for linear mappings.
A semiring $K$ or the normal completion $\wh K$
of a semifield $K$ is a semimodule over $K$. If $W$
coincides with $K$ (or $\wh K$), then the linear mapping $p$
is called a {\it linear functional}. Linear
(respectively, $b$-linear) mappings will also be referred to as
{\it linear\/} (respectively, $b$-{\it linear\/})
operators.

In analysis, the most important examples of idempotent
semimodules and spaces are either subsemimodules of
topological vector lattices~\cite{15}
(possibly coinciding with these) or their duals (that is,
semimodules of linear
functionals with some regularity condition, say, of $a$-linear
functionals).

\remark{Remark~1}
Clearly, idempotent semimodules over a given idempotent semiring
form a category with linear mappings as
morphisms. In the present paper we mainly deal with another
category, which is formed by $b$-complete idempotent
semimodules over a given $b$-complete idempotent semiring with
$b$-linear mappings as morphisms, and also with the
full subcategory of $a$-complete
semimodules with $a$-linear mappings as
morphisms. Note that these categories are not additive
(it is natural to say that they are {\it semiadditive\/});
for the basic notions of the category theory, e.g., see~\cite{16}.
\endremark

In what follows, unless specified otherwise, we consider only
$b$-complete (including $a$-complete)
idempotent semirings and idempotent semimodules.

\head
\S2. Direct Products and Sums of  $b$-Complete
Idempotent Semimodules
\endhead

Let $K$ be a $b$-complete idempotent semiring,
and let $\{V_\alpha\}_{\alpha\in A}$ be a family of
$b$-complete
idempotent semimodules over $K$. We represent an element
of the direct
product of all sets of this family as a (generally speaking,
infinite) sequence (``vector'')
$x = \{x_\alpha\} = (\dots,x_\alpha,\dots)$,
where
$x_\alpha\in V_\alpha$. The componentwise operations
$x\oplus y = \{x_\alpha+y_\alpha\} =
(\dots,x_\alpha\oplus y_\alpha,\dots)$ and
$k\odot x = \{k\odot x_\alpha\} =
(\dots,k\odot x_\alpha,\dots)$,
$k\in K$, make the direct product of the sets
$V_\alpha$, $\alpha\in A$, an idempotent semimodule,
which is called the {\it direct product
of the semimodules\/} $V_\alpha$, $\alpha\in A$, and is
denoted by $\prod_\alpha V_\alpha$. The direct product
of semimodules $V_1,\dots,V_n$ will also be denoted by
$V_1\times\dots\times V_n$. For each index
$\alpha\in A$, we define a projection
$p_\alpha\:\prod_\alpha V_\alpha\to V_\alpha$ by setting
$p_\alpha(x) = x_\alpha$ and an embedding
$i_\alpha\:V_\alpha\to\prod_\alpha V_\alpha$ by setting
$(i(x))_\alpha = x$ and $(i(x))_\beta = \0$ for $\alpha\ne\beta$.

\proclaim{Proposition~1}
The semimodule $\prod_\alpha V_\alpha$ is $a$-complete
\rom(respectively, $b$-complete\rom) if so are all semimodules
$V_\alpha$.
\endproclaim

This assertion readily follows from the definitions.

\proclaim{Proposition~2}
Let $\{V_\alpha\}_{\alpha\in A}$ be a family of $b$-complete
semimodules over $K$, and let $V$ be a $b$-complete semimodule
over $K$. Then the following assertions hold\rom:
\roster
\item"\rom{1)}" for each family of $b$-linear mappings
$f_\alpha\:V\to V_\alpha$ there exists a unique
$b$-linear mapping $f\:V\to\prod_\alpha V_\alpha$
such that $p_\alpha f = f_\alpha$ for each index
$\alpha\in A$\rom;
\item"\rom{2)}" if the indexing set $A$ is finite, then for
each family of $b$-linear mappings
$f_\alpha\:V_\alpha\to V$ there exists a unique
$b$-linear mapping $f\:\prod_\alpha V_\alpha\to V$
such that $fi_\alpha = f_\alpha$ for each index
$\alpha\in A$\rom;
\item"\rom{3)}" if all semimodules $V$ and $V_\alpha$,
$\alpha\in A$, are $a$-complete, then for any
indexing set $A$ and any family of $a$-linear
mappings $f_\alpha\:V_\alpha\to V$ there exists
a unique $a$-linear mapping
$f\:\prod_\alpha V_\alpha\to V$ such that
$fi_\alpha = f_\alpha$ for each $\alpha\in A$.
\endroster
\endproclaim

\demo{Proof}
We start from assertion~1). We set $f(x) = \{f_\alpha(x)\}$.
One can readily see that a subset
of $\prod_\alpha V_\alpha$ is bounded if and only if
the projections of this subset on all
components are bounded. This, together with the definitions,
readily implies that $f$ is $b$-linear; moreover, by
construction,
$p_\alpha f = f_\alpha$ for all $\alpha\in A$.
The uniqueness of the desired mapping is also obvious, and the
proof of assertion~1) is complete. The mapping $f$ is called the
{\it direct product\/} of the mappings $f_\alpha$.

To prove assertion~2), we set
$f(\{x_\alpha\}) = \bigoplus_\alpha f_\alpha(x_\alpha)$.
One can readily see that if $X\subset \prod_\alpha V_\alpha$
is a bounded subset, then so is the set
$\{f_\alpha(x_\alpha)\mid\alpha\in A,\ x_\alpha\in p_\alpha(X)\}$.
It follows that the mapping $f$ is well defined; the relations
$fi_\alpha = f_\alpha$ for all $\alpha\in A$ and the fact that
$f$ is $b$-linear and unique can be verified directly.
This mapping is called the {\it direct sum\/} of the
mappings $f_\alpha$.

The proof of assertion~3) is similar. The finiteness of
the indexing set was used only in the proof of the
boundedness of the set
$\{f_\alpha(x_\alpha)\mid\alpha\in A,\
x_\alpha\in p_\alpha(X)\}$.
However, for an $a$-complete semimodule $V$ this holds
automatically. It remains to notice that under the assumptions
of this
assertion $a$-linearity coincides with $b$-linearity.

Thus the proof of Proposition~2 is complete.
\qed\enddemo

\remark{Remark~2}
In the language of the category theory, Proposition~2 has the
following meaning. In the
category of $b$-complete idempotent semimodules there exists a
(categorical) direct product of an arbitrary family of
semimodules $V_\alpha$; this product coincides with
$\prod_\alpha V_\alpha$.  If the family is finite, then
there also exists a (categorical) direct sum
$\sum_\alpha V_\alpha$ semimodules $V_\alpha$, which coincides
with the direct sum $\prod_\alpha V_\alpha$. In the category
of $a$-complete semimodules, for an arbitrary family
of semimodules $V_\alpha$ there exist a (categorical) direct
sum $\sum_\alpha V_\alpha$ and a direct
product $\prod_\alpha V_\alpha$, which coincide with each other and
with the direct product in the category of $b$-complete
semimodules.
\endremark

\head
\S3. Tensor Products of $b$-Complete Semimodules
\endhead

Let $\{V_\alpha\}_{\alpha\in A}$ be a family of $b$-complete
idempotent semimodules over a given $b$-complete
commutative idempotent semiring $K$, and let
$V = \prod_\alpha V_\alpha$.

By $T$ we denote the set of formal sums of the form
$$
t = \bigoplus_{x = \{x_\alpha\}\in X}\lambda(x)
\odot\bigotimes_\alpha x_\alpha,
\tag1
$$
where $X\subset V$, $x = \{x_\alpha\} = (\dots,x_\alpha,\dots)$,
and $\lambda$ is an arbitrary mapping of $X$
into $K$. An element of the form~\thetag{1} will be called a
{\it representation\/} of a tensor (in what follows we define a
tensor as a class of equivalent representations).
The natural (formal) idempotent addition (which should not be
confused with the addition in $V = \prod_\alpha V_\alpha$) and
the multiplication by elements $k\in K$ (which takes each
function $\lambda(x)$ to $x\mapsto k\odot\lambda(x)$) make $T$
a semimodule over $K$.

We say that a representation of the form~\thetag{1} is {\it
bounded\/} if $X$ is bounded in $V$ and $\lambda(x)$
is a bounded function on $X$ (ranging in $K$).
The set of all bounded representations will be denoted by
$T_b$. Clearly, $T_b$ is a subsemimodule of $T$. One can
readily see that the semimodule $T_b$ is
$b$-complete.

As usual, we equip $T$ and $T_b$ with the equivalence
relation generated by the identities
$$
\align
k\odot(\dotsb\otimes x_\alpha\otimes\dotsb) =
& (\dotsb\otimes k\odot x_\alpha\otimes\dotsb),
\tag2
\\
(\dotsb\otimes(\oplus X_\alpha)\otimes\dotsb) =
& \bigoplus_{x\in X_\alpha}
(\dotsb\otimes x\otimes\dotsb),
\tag3
\endalign
$$
where $k\in K$, $X_\alpha\subset V_\alpha$, some value of the
index $\alpha$ is chosen, and all respective components
(formal factors) on the right- and left-hand sides
in~\thetag{2} and~\thetag{3} coincide except for
those written out explicitly. Relation~\thetag{2} permits one
to replace the sum~\thetag{1} by an equivalent
representation of the form
$$
t = \bigoplus_{x = \{x_\alpha\}\in X}\bigotimes_\alpha x_\alpha.
\tag4
$$

In what follows we deal with representations
of the form~\thetag{4}; in this case, by some abuse of speech,
the set $X\subset\prod_\alpha V_\alpha$ occurring
in~\thetag{4} will also be called a {\it representation\/} of
a given tensor. {\it The addition of tensors of the form\/}~\thetag{4}
{\it corresponds to the union
of their representations \rom(summands\rom)}. This addition
must be idempotent; hence, by summing a given
representation with all equivalent representations, we obtain
the same tensor (up to equivalence).
{\it The corresponding representation is
naturally called complete\rom; we identify a tensor\rom,
i.e., a class of equivalent representations\rom, with the
corresponding complete
representation}. We shall assume that the complete representation
always contains the element
$\0_V = \bigotimes_\alpha\0_{V_\alpha}\in V$.

Let us proceed to precise definitions. For each element
$k\in K$, we denote by $k_\alpha$ the self-mapping of
$V = \prod_\alpha V_\alpha$ that acts as the multiplication by
$k$ of the $\alpha$\snug th component of each element $x\in V$
and does not alter any other components.
A subset $S\subset V$ will be called an $\alpha$-{\it fiber\/} or
simply a {\it fiber\/} if the projection $p_\beta(S)$ is a
singleton for $\alpha\ne\beta$ and if $p_\alpha(S) = V_\alpha$.
Thus, all but one components of elements of a fiber $S$ are
fixed, whereas the $\alpha$\snug th component can take
arbitrary values in $V_\alpha$.

A {\it tensor}, or a {\it complete representation\/} of a
tensor, is an arbitrary subset $X\subset V =
\prod_\alpha V_\alpha$ with the following properties:
\roster
\item"1)" if $k_\alpha(x)\in X$ for some index
$\alpha\in A$ and $k\in K$, then $k_\beta(x)\in X$ for all
indices $\beta\in A$;
\item"2)" $\oplus(X\cap S)\in X$ for each fiber
$S\subset V$;
\item"3)" if $y\in S$, $x\in X\cap S$, and
$y\preccurlyeq x$, then $y\in X$ for each fiber $S\subset V$.
\endroster

Needless to say, property~1) corresponds to
Eq.~\thetag{2}, and properties~2) and~3) correspond to
Eq.~\thetag{3} with regard to the fact that
$x\oplus y = x$ for
$y\preccurlyeq x$. We shall assume that a complete
representation $X$ always contains $\0_V = \{\0_\alpha\}$.
Formally, this corresponds to the identity $\sup\emptyset = \0_V$.

We define the $\tau$-{\it hull\/} $X^\tau$
of an arbitrary subset $X\subset V = \prod_\alpha V_\alpha$
as the intersection of all
tensors (that is, their complete representations) containing
$X$. Two subsets $X,Y\subset V$ are said to be
{\it equivalent\/} (we denote this equivalence
relation by $\tau$) if their $\tau$-hulls
$X^\tau$ and $Y^\tau$ coincide.

An elementary analysis shows that the following assertion holds.

\proclaim{Proposition~3}
The equivalence generated by~\thetag{2}
and~\thetag{3} with regard to the fact that the addition is idempotent
on the set of representations of the form~\thetag{4}
coincides with the equivalence $\tau$. Every representation
of the form~\thetag{1} is equivalent to a unique complete
representation of the form~\thetag{4}. This equivalence
is consistent with the structure of a semimodule over $K$ on
$T$ and $T_b$.
\endproclaim

It follows that the quotient semimodules $\wt T$
and $\wt T_b$ are well defined. One can readily see that
the semimodule $\wt T_b$
is  $b$-complete; in any case, this follows from
Proposition~4 below. The quotient semimodule $\wt T_b$ will be denoted
by $\bigotimes_\alpha V_\alpha$ and called the
$b$-{\it tensor\/} (or simply {\it tensor\/})
{\it product of the $b$-complete semimodules\/} $V_\alpha$.

Recall that we identify a tensor with its complete
representation, i.e., a subset $X X^\tau\subset V$.

We say that a tensor $X$ is a {\it bounded\/} or
$b$-{\it tensor} if it is contained (as a set) in a
tensor $\{x\}^\tau$, where $x\in V$. In this case we say that
$X$ is bounded by the tensor $\{x\}^\tau$, i.e.,
the tensor with representation $\bigotimes_\alpha x_\alpha$,
where $x = \{x_\alpha\}$.

Let $T_b(V)$ be the set of all $b$-tensors,
ordered by inclusion. This order gives rise to an
idempotent addition in $T_b(V)$. The sum of an arbitrary
bounded family $\{X_\alpha\}$ of tensors coincides with
the $\tau$-hull $\{\bigcup_\alpha X_\alpha\}^\tau$ of their
unions. If a $b$-tensor $X$ is bounded
by a tensor $\{x\}^\tau$ and a $b$-tensor $Y$ is bounded
by a tensor $\{y\}^\tau$, then the tensor
$X\oplus Y = \{X\cup Y\}^\tau$ is bounded by the tensor
$\{x\oplus y\}^\tau$, where $x\oplus y\in V$. Clearly, the
intersection of an arbitrary family of $b$-tensors is again a
$b$-tensor, which coincides with the greatest lower bound of this
family of tensors (with respect to the order defined by the
set-theoretic inclusion of subsets of $V$). It follows
that every subset $M\subset T_b(V)$ that is bounded above
has the least upper bound
$$
\sup M = \inf\{X\mid Y\preccurlyeq X\ \text{for all}\ Y\in M\}
= \bigcap_{X\supset Y\in M}X
$$
(e.g., see~\cite{10}).
Thus $T_b(V)$
is a lattice and a $b$-complete idempotent semigroup.

One can multiply elements $X\in T_b(V)$ by elements
 $k\in K$ according to the formula
$k\odot X = (k_\alpha(X))^\tau$ for
a given index $\alpha\in A$; by property~1)
of the complete representation of a tensor, this
multiplication is independent of the choice of $\alpha$.

\proclaim{Proposition~4}
The idempotent semigroup $T_b(V)$ with the above-defined
multiplication by elements of $K$ is a $b$-complete
semimodule over $K$. This semimodule is isomorphic to the
semimodule
$\bigotimes_\alpha V_\alpha$.
\endproclaim

\demo{Proof}
Let us verify that $T_b(V)$ is a $b$-complete semimodule.
We have already seen that $T_b(V)$ is a $b$-complete idempotent
semigroup.
Hence it suffices to verify the distributivity
laws and the associativity of the multiplication by elements of $K$.
One can readily verify that for each $b$-homomorphism
$g\:V\to V$ and each $b$-tensor $X\subset V$,
the set $g^{-1}(X)$ is also a $b$-tensor. It follows that
$g(X^\tau)\subset(g(X))^\tau$ and
$(g(X^\tau))^\tau = (g(X))^\tau$
for each subset $X\subset V$. By applying these relations to
$b$-homomorphisms of the form $k_\alpha$, we find that
$(k_\alpha(X))^\tau = (k_\alpha(X^\tau))^\tau$ for all
$X\subset V$, $k\in K$, $\alpha\in A$, whence it follows
that the multiplication by elements of $K$ is associative.

Let us verify the distributivity. For each subset $\sT\subset
T_b(V)$ that is bounded above, one has
$\oplus\sT = \bigl(\,\bigcup_{t\in\sT}t\bigr)^\tau$.
Consequently,
$$
k\odot(\oplus\sT)=
 k\odot\biggl(\,\bigcup_{t\in\sT}t\biggr)^\tau=
 \biggl(k_\alpha\biggl(\biggl(\,\bigcup_{t\in\sT}t
\biggr)^\tau\biggr)\biggr)^\tau=
 \biggl(k_\alpha\biggl(\,\bigcup_{t\in\sT}t\biggr)\biggr)^\tau=
 \biggl(\,\bigcup_{t\in\sT}k_\alpha(t)\biggr)^\tau.
$$
Since
$$
\biggl(\,\bigcup_{t\in\sT}k_\alpha(t)\biggr)^\tau=
 \biggl(\,\bigcup_{t\in\sT}(k_\alpha(t))^\tau\biggr)^\tau=
 \bigoplus_{t\in\sT}k\odot t = \oplus(k\odot\sT),
$$
we find that $k\odot(\oplus\sT) = \oplus(k\odot\sT)$, so
that one of the distributivity laws follows. Here the small
letter $t$ stands for a $b$-tensor, whereas above (and below)
tensors are denoted by capital letters.

Let us now verify the second distributivity law
$(\oplus Q)\odot X = \oplus(Q\odot X)$, where $Q$ is an
arbitrary bounded subset in $K$ and $X$ is an arbitrary
$b$-tensor. It suffices to verify that the
intersections of these sets with the $\alpha$-fiber $S$
coincide for
an arbitrary index $\alpha\in A$. However, by construction,
$(\oplus Q)_\alpha(x) = \oplus(Q\odot x)\in\{k\odot x\mid k\in Q\}^\tau$
for any $x\in S$.
Consequently, $(\oplus Q)_\alpha(x)\in
\bigl(\,\bigcup_{k\in Q}k_\alpha(X)\bigr)^\tau$
for $x\in S\cap X$, whence it follows that
$((\oplus Q)_\alpha(X))\subset
\bigl(\,\bigcup_{k\in Q}k_\alpha(X)\bigr)^\tau$
and hence
$(\oplus Q)\odot X = ((\oplus Q)_\alpha(X))^\tau\subset
\bigl(\,\bigcup_{k\in Q}k_\alpha(X)\bigr)^\tau$.
However,
$\bigl(\,\bigcup_{k\in Q}k_\alpha(X)\bigr)^\tau =
\bigl(\,\bigcup_{k\in Q}(k_\alpha(X))^\tau\bigr)^\tau =
\bigoplus_{k\in Q}k\odot X = \oplus(Q\odot X)$.
Thus
$(\oplus Q)\odot X\preccurlyeq\oplus(Q\odot X)$. Since
the opposite inequality is obvious, this completes the
verification of the second distributivity law. Thus we have
proved that $T_b(V)$
is a $b$-complete semimodule over $K$.

It remains to indicate (and verify) an isomorphism between
$T_b(X)$ and $\bigoplus_\alpha V_\alpha$. Basically, now we
can readily derive this from Proposition~3. We obtain a mapping
$T_b(V)\to\bigoplus_\alpha V_\alpha$ by taking each
$b$-tensor $X$ to its (complete) formal representation
$\bigoplus_{x = \{x_\alpha\}\in X}\bigotimes_\alpha x_\alpha$
of the form~\thetag{4} and by passing to the corresponding
equivalence class in $T_b$. The inverse mapping
takes an arbitrary representation of the form~\thetag{1}
or~\thetag{4} to the complete representation of the form~\thetag{4}.
The verification of the fact that these mappings are mutually
inverse and the structures of semimodules over $K$ are
consistent (coincide) is elementary. The proof of the
proposition is complete.
\qed\enddemo

\remark{Remark~3}
Relation~\thetag{2} and property~1) of the complete representation of a
tensor show why we require the basic semiring $K$ to be commutative.
\endremark

\head
\S4. Basic Results on $b$-Polylinear Mappings
\endhead

Let $\{V_\alpha\}_{\alpha\in A}$ be an arbitrary
family of $b$-complete idempotent semimodules over a
$b$-complete commutative semiring $K$, and let $W$ be
an arbitrary $b$-complete semimodule over $K$. We set
$V = \prod_\alpha V_\alpha$; thus, $V$ is the direct
product of the semimodules $V_\alpha$; by
$\bigotimes_\alpha V_\alpha$ we denote the $b$-tensor
product of these semimodules.

A mapping $f\:\prod_\alpha V_\alpha\to W$ is said to be
$b$-{\it polylinear\/} if it is separately $b$-linear in
every component $V_\alpha$; this means that for every
given index $\beta\in A$ the mapping
$x_\beta\mapsto f(\dots,x_\beta,\dots)$ for
arbitrary fixed values of all other coordinates is a
$b$-linear mapping $V_\beta\to W$.

We define a {\it canonical mapping\/} $\pi$
{\it of the direct product\/} $V = \prod_\alpha V_\alpha$
{\it into the  tensor product\/}
$T_b(V) = \bigotimes_\alpha V_\alpha$ by setting
$\pi(x) = \{x\}^\tau$, i.e.,
$\pi(\{x_\alpha\}) = \bigotimes_\alpha x_\alpha$.

\remark{Remark~4}
The range of the canonical mapping $\pi$ generates the
semimodule
$\bigotimes_\alpha V_\alpha$, i.e., each element of this
semimodule is a linear combination (not necessarily
finite) of elements of the form $\bigotimes_\alpha x_\alpha$;
moreover, each element is a sum (not necessarily
finite) of elements of the form $\bigotimes_\alpha x_\alpha$.
In other words, the set $\pi(\prod_\alpha V_\alpha)$
is a {\it system of generators\/} of the semimodule
$\bigotimes_\alpha V_\alpha$ and even of the
corresponding idempotent semigroup.
\endremark

\proclaim{Theorem~1}
The canonical mapping
$\{x_\alpha\}\mapsto\bigotimes_\alpha x_\alpha$ of the direct
product $V = \prod_\alpha V_\alpha$
of the semimodules $V_\alpha$ into the  tensor product
$T_b(V) = \bigotimes_\alpha V_\alpha$ of these semimodules
is $b$-polylinear. For each $b$-polylinear
mapping $f\:\prod_\alpha V_\alpha\to W$ there exists
a unique $b$-linear mapping
$f_\otimes\:\bigotimes_\alpha V_\alpha\to W$ such that
$f = f_\otimes\pi$.
\endproclaim

To prove this theorem, we need some
auxiliary assertions. Just as above, we identify a
tensor with its complete representation
of the form~\thetag{4}, i.e., with the corresponding subset
in $\prod_\alpha V_\alpha$. For each element $w$ of
an arbitrary idempotent semigroup $W$, we denote by
$\Low(w)$ the set
$\{u\in W\mid u\preccurlyeq w\}$.

\proclaim{Lemma~1}
Let $W$ be an arbitrary $b$-complete semimodule over $K$.
For each $b$-polylinear mapping
$f\:\prod_\alpha V_\alpha\to W$, the complete preimage of
each
set of the form $\Low(w)$, where $w\in W$, is a
$b$-tensor.
\endproclaim

\demo{Proof}
Let $X = f^{-1}(\Low(w))$. If $x$ and $y$ are elements
of some fiber, $x\in X$, and $y\preccurlyeq x$, then
$f(y)\preccurlyeq f(x)\preccurlyeq w$, whence $y\in X$.
If $S\subset X$ and $S$ is contained in a single fiber, then
$f(\oplus S) = \oplus f(S)\preccurlyeq w$, whence
$\oplus S\in X$. If $y = k_\alpha(x)$ and $k_\beta(x)\in X$
for some index $\beta$, then
$f(y) = k\odot f(x) = f((k_\beta(x))\preccurlyeq w$, whence
$y\in X$. Thus conditions~1)--3) from the definition of a
tensor are satisfied, which completes the proof of Lemma~1.
\qed\enddemo

\proclaim{Lemma~2}
For any $b$-polylinear mapping
$f\:\prod_\alpha V_\alpha\to W$ and any subset $X\subset
V = \prod_\alpha V_\alpha$, one has
$\oplus f(X) = \oplus f(X^\tau)$.
\endproclaim

\demo{Proof}
It follows from Lemma~1 that the set
$\{x\in V\mid f(x)\preccurlyeq w\}$ is a $b$-tensor
for each $w\in W$. In particular, the set
$\{x\in V\mid f(x)\preccurlyeq\oplus f(X)\}$ is a
$b$-tensor containing $X$. It follows that it contains
the set $X^\tau$, whence
$\oplus f(X)\succcurlyeq\oplus f(X^\tau)$. Since
the opposite inequality is obvious, the proof of Lemma~2 is
complete.
\qed\enddemo

\demo{Proof of Theorem 1}
Let us verify that the canonical mapping is $b$-polylinear.
To this end, let us verify that
$\pi(\oplus X) = \oplus\pi(X)$
for each subset $X\subset
V = \prod_\alpha V_\alpha$ lying in an $\alpha$-fiber
$S\subset V$.
Since
$\oplus X\in X^\tau$, we have the inclusion
$\pi(\oplus X)\subset X^\tau$. On the other hand, if
$x\in X$, then $x\in S$, $\oplus X\in S$, and
$x\preccurlyeq\oplus X$, whence $x\in(\oplus X)^\tau$.
Consequently, $X\subset(\oplus X)^\tau$, whence
$X^\tau\subset(\oplus X)^\tau = \pi(\oplus X)$, so that
$\oplus\pi(X) = X^\tau = \pi(\oplus X)$, as desired.
Since the relation $\pi(k_\alpha(x)) = k\odot\pi(x)$, where
$x\in V$, $k\in K$, readily follows from
the definition of the multiplication of a $b$-tensor by an
element $k\in K$,
we have proved that the canonical mapping $\pi$ is $b$-polylinear.

It remains to verify the existence and uniqueness of
the mapping $f_\otimes$. For each $b$-tensor $t$,
we set $f_\otimes(t) = \oplus f(t)$, where $f(t)$ is the image
of the set $t$ (the complete representation of the tensor $t$)
under the
mapping $f$. It follows from Lemma~2 that this is well-defined
and that $f_\oplus(\pi(x)) = f(x)$.

For each subset $T$ of the tensor
product $\bigotimes V_\alpha$, one has
$$
f_\otimes(\oplus T)=
 f_\otimes\biggl(\biggl(\,\bigcup_{t\in T}t\biggr)^\tau\biggr)=
 \oplus f\biggl(\,\bigcup_{t\in T}t\biggr)=
 \bigoplus_{t\in T}(\oplus f(t))=
 \oplus f_\otimes(T),
$$
so that $f_\otimes$ is a $b$-homomorphism of the
corresponding idempotent semigroups.

One can also readily prove that the mapping $f_\otimes$ is
homogeneous. Indeed,
$$
f_\otimes(k\odot\pi(x))=
 f_\otimes(\pi(k_\alpha(x)))=
 f(k_\alpha(x))=
 k\odot f(x)=
 k\odot f_\otimes(\pi(x)).
$$
Thus $f_\otimes$ is homogeneous on $\pi(V)$.
Since $f_\otimes$ is a $b$-homomorphism and
the range $\pi(V)$ of the canonical mapping
generates $\bigotimes_\alpha V_\alpha$ (see Remark~4),
this completes the proof of the homogeneity of the mapping
$f_\otimes$.
The uniqueness of this mapping also follows from
Remark~4. The proof of Theorem~1 is complete.
\qed\enddemo

\proclaim{Theorem~2}
Suppose that $V = \prod_\alpha V_\alpha$, $W$ is a $b$-complete
semimodule over $K$, and there is a $b$-polylinear
mapping $\epsilon\:V\to W$ whose range
generates the semimodule $W$, so that each element of $W$
is a linear combination of elements
of the form $\epsilon(x)$. Suppose that for each
$b$-polylinear mapping $f\:V\to U$, where $U$
is an arbitrary $b$-complete semimodule over $K$, there exists
a $b$-linear mapping $F\:W\to U$ such that
$F\epsilon = f$. Then there exists an isomorphism
$\delta\:W\to\bigotimes_\alpha V_\alpha$ of the semimodule $W$
onto the
$b$-tensor product $\bigotimes_\alpha V_\alpha$
such that $\delta\epsilon = \pi$.
\endproclaim

\demo{Proof}
Let $\delta\:W\to\bigotimes_\alpha V_\alpha$ be a
$b$-linear mapping such that $\delta\epsilon = \pi$.
This mapping exists by the assumption of Theorem~2 applied to
the $b$-polylinear mapping $f = \pi$. On the other hand, by
applying Theorem~1 to the $b$-polylinear
mapping $\epsilon$, we can construct a $b$-linear mapping
$\epsilon_\otimes\:\bigotimes_\alpha V_\alpha\to W$ such
that $\epsilon_\otimes\pi = \epsilon$. Since
$\delta\epsilon_\otimes\pi = \delta\epsilon = \pi$, it follows
that the
$b$-linear mapping $\delta\epsilon_\otimes$
leaves all elements of $\pi(V)$ in their places; now we see
that this
mapping coincides with the identity mapping, since $\pi(V)$
generates $\bigotimes_\alpha V_\alpha$. Likewise, one can
check that $\epsilon_\otimes\delta$ is also the
identity mapping. Thus, the $b$-linear mappings
$\epsilon_\otimes$ and $\delta$ are the inverses of each other
and are isomorphisms. Thus Theorem~2 follows from Theorem~1.
\qed\enddemo

\remark{Remark~5}
Theorems~1 and~2 show that in the category of $b$-complete
semimodules (with $b$-linear mappings as
morphisms) we have constructed a natural tensor product.
This tensor product is also a natural tensor product in the
full subcategory of $a$-complete
semimodules, since one can readily see that the $b$-tensor
product of an arbitrary family of $a$-complete semimodules (over a
given idempotent semiring) is an $a$-complete semimodule.
\endremark

By using similar constructions,
one can readily construct a tensor product (with finitely many
factors) in the category of arbitrary idempotent
semimodules over a given commutative idempotent semiring
(with linear mappings as morphisms).

\head
\S5. Tensor Products of $b$-Linear Mappings
\endhead

Let $\{V_\alpha\}_{\alpha\in A}$ and
$\{W_\alpha\}_{\alpha\in A}$ be families of $b$-complete
semimodules over a given $b$-complete commutative
semiring $K$. Suppose that for each index $\alpha\in A$
there is a $b$-linear mapping
$f_\alpha\:V_\alpha\to W_\alpha$. By $f$ we denote the
direct product $\prod_\alpha f_\alpha$ of these
mappings, i.e., the mapping
$f\:\prod_\alpha V_\alpha\to\prod_\alpha W_\alpha$ such that
$f(\{x_\alpha\}) = \{f(x_\alpha)\}$. One can readily see that this
mapping extends to be a $b$-polylinear mapping
$$
\prod_\alpha V_\alpha\to\bigotimes_\alpha W_\alpha.
\tag5
$$
It follows from Theorem~1 that the mapping~\thetag{5}
extends to be a $b$-linear mapping
$$
\bigotimes V_\alpha\to\bigotimes_\alpha W_\alpha.
\tag6
$$
The mapping~\thetag{6} is called the $b$-{\it tensor\/}
(or simply {\it tensor\/}) {\it product\/}
of the $b$-{\it linear mappings\/} $f_\alpha$ and is
denoted by $\bigotimes_\alpha f_\alpha$.

\head
\S6. The Tensor Algebra
\endhead

Let $U,V$ and $W$ be $b$-complete semimodules over a $b$-complete
commutative idempotent semiring $K$. We use the symbol
$\dotplus$ to denote the direct sum of two
semimodules and the symbol $\sum_\alpha V_\alpha$ to
denote the direct sum of an arbitrary family
of $a$-complete semimodules. In these cases the direct sum
coincides with the direct product (see Remark~2).

\proclaim{Theorem~3}
The following $b$-complete semimodules are isomorphic\rom:
\roster
\item"\rom{1)}" $U\otimes V$ and $V\otimes U$\rom;
\item"\rom{2)}" $(U\dotplus V)\otimes W$ and
$U\otimes W\dotplus V\otimes W$\rom;
\item"\rom{3)}" $(U\otimes V)\otimes W$ and $U\otimes(V\otimes W)$\rom;
\item"\rom{4)}" $(\,\sum_\alpha V_\alpha)\otimes W$ and
$\sum_\alpha V_\alpha\otimes W$,
\endroster
where $\{V_\alpha\}_{\alpha\in A}$ is an arbitrary family of
$a$-complete semimodules over $K$.
\endproclaim

\proclaim{Corollary}
The set of isomorphism classes of $b$-complete semimodules \rom(over a
given $b$-complete commutative~\rom{idempotent semiring}\rom)
is a commutative
associative semiring with respect to the operations  of direct
sum and $b$-tensor product.
\endproclaim

\remark{Remark~6}
This semiring of semimodules is not idempotent.
\endremark

\demo{Proof of Theorem 3}
Assertion~1) is trivial and readily follows from the
definition of the $b$-tensor product. Assertion~3)
follows from the fact that the direct products
$U\times(V\times W)$ and $(U\times V)\times W$ are isomorphic,
and moreover, the set of $b$-tensors (i.e., of their complete
representations of the form~\thetag{4}) is the same in both
cases and
coincides with the set of $b$-tensors in $U\times V\times W$,
i.e., with $U\otimes V\otimes W$, as desired.

Let us prove assertions~2) and~4). Suppose that the direct sum
$\sum_\alpha V_\alpha = V$ of the family of semimodules
$V_\alpha$ exists and coincides with the direct
product $\prod_\alpha V_\alpha$ (see Remark~2);
for assertion~2) this family is finite (or even consists
of two semimodules). Let $i_\alpha\:V_\alpha\to V$ and
$p_\alpha\:V\to V_\alpha$ be the corresponding  canonical
embeddings and projections, and let
$i\:\sum_\alpha(V_\alpha\otimes W)\to V\otimes W$ be the
direct sum of the mappings $i_\alpha\otimes I_W$, where
$I_W$ is the identity operator in $W$ (the tensor
product of $b$-linear mappings was defined in~\S5).
Thus the restriction of the mappings $i$
to $V_\alpha\otimes W$ coincides with $i_\alpha\otimes I_W$.
Likewise, let
$p\:V\otimes W\to\sum_\alpha(V_\alpha\otimes W) =
\prod_\alpha(V_\alpha\otimes W)$
be the direct product of the
mappings $p_\alpha\otimes I_W$. A routine verification shows
that $i$ and $p$ are mutually inverse isomorphisms.
\qed\enddemo

\head
\S7. Tensor Products of Semimodules of Bounded Functions
\endhead

Let $K$ be a $b$-complete commutative idempotent semiring and
$X$ an arbitrary set. By $\sB(X,K)$ we denote
the set of {\it bounded functions on\/} $X$ {\it ranging in\/}
$K$, i.e., mappings $X\to K$ with
bounded range. If $f,g\in\sB(X,K)$,
then the pointwise addition and multiplication
$$
(f\oplus g)(x) = f(x)\oplus g(x),
\qquad
(f\odot g)(x) = f(x)\odot g(x),
$$
where $x$ is an arbitrary element of $X$, define a structure
of a $b$-complete semimodule over $K$ on $\sB(X,K)$
(e.g., see~\cite{17, 18, 5, 12}).

\proclaim{Proposition~5}
For any sets $X$ and $Y$, the $b$-complete semimodules
$\sB(X\times Y,K)$
and $\sB(X,K)\otimes\sB(Y,K)$ are isomorphic.
\endproclaim

\demo{Proof}
Let $\delta_x\in\sB(X,K)$ be the function equal to
$\1$ at the point $x\in X$ and $\0$ at all other points. One
can readily verify that the mapping
$$
(f,g)\mapsto\oplus\{f(x)\odot f(y)\odot\delta_{(x,y)}\}
$$
is a $b$-linear mapping
$\sB(X,K)\times\sB(Y,K)\to\sB(X\times Y,K)$;
hence it extends to be a $b$-linear mapping
$$
i\:\sB(X,K)\otimes\sB(Y,K)\to\sB(X\times Y,K).
\tag7
$$

We define a $b$-linear mapping
$$
j\:\sB(X\times Y,K)\to\sB(X,K)\otimes\sB(Y,K)
\tag8
$$
by setting
$j(f) = \oplus\{f(x,y)\odot\delta_x\otimes\delta_y\}$; this is
obviously well defined.
One can readily see that the mappings~\thetag{7} and~\thetag{8}
are the inverses of each other and specify the desired homomorphism.
It suffices to verify this on a system of generators. For an
element $(x,y)\in X\times Y$, one readily has
$i(j(\delta_{(x,y)})) = i(\delta_x\otimes\delta_y) = \delta_{(x,y)}$,
so that $ij$ is the identity mapping
on $\sB(X\times Y,K)$. On the other hand, elements of the form
$\delta_x\otimes\delta_y$ obviously generate
$\sB(X,K)\otimes\sB(Y,K)$, and
$j(i(\delta_x\otimes\delta_y)) =
j(\delta_{(x,y)}) = \delta_x\otimes\delta_y$;
it follows that $ji$ is the identity mapping.
The proof of the proposition is complete.
\qed\enddemo

\head
\S8. Nuclear Operators on $b$-Complete Semimodules
\endhead

Let $V$ and $W$ be $b$-complete semimodules over  a
$b$-complete commutative idempotent semiring $K$. By $W^*$ we
denote the
semimodule of all $b$-linear mappings $W\to K$
($b$-linear functionals). The set of all $b$-linear
mappings $W\to V$ is denoted by $\Hom_b(W,V)$.
This set bears a natural structure of a $b$-complete
semimodule over $K$. Thus $W^* = \Hom_b(W,K)$.

For $v\in V$ and $w^*\in W^*$ we define an operator (i.e., a
$b$-linear mapping) $p_{w^*,v}\: W\to V$
by the formula
$$
p_{w^*,v}\:x\mapsto w^*(x)\odot v.
\tag9
$$
Operators of the form~\thetag{9} will be called {\it
operators of rank\/}~1 (or {\it one-dimensional operators\/}).

\proclaim{Proposition~6}
There exists a unique $b$-linear mapping
$p\:W^*\otimes V\to\Hom_b(W,V)$ such that
$p(w^*\otimes v) = p_{w^*,v}$. The range of the operator $p$ is
the
subsemimodule of $\Hom_b(W,V)$ generated by the set of
all operators of the form~\thetag{9}, i.e., by operators of
rank~\rom1.
\endproclaim

To prove this, it suffices to notice that
the mapping
$W^*\times V\to\Hom_b(W,V)$ taking each pair $(w^*,v)$ to the
operator $p_{w^*,v}$ is $b$-bilinear and its range
consists of all operators of rank~1; then we can apply
Theorem~1. The mapping $p$ described in Proposition~6 will be
called {\it canonical}.

By analogy with~\cite{1},
we refer to a $b$-linear mapping
$n\:W\to V$ as a $b$-{\it nuclear operator\/} (or
a {\it nuclear mapping\/}) if $n$ lies in the range
of the canonical mapping $p$. It follows from Proposition~6
that a nuclear mapping $W\to V$ can be represented as a
sum of operators of rank~1. It readily follows from this
proposition and from our definitions that the following
assertion holds.

\proclaim{Corollary}
If the semimodules $V$ and $W$ are $a$-complete, then
the range
of the
canonical mapping $p$ is an $a$-complete
subsemimodule of $\Hom_b(W,V)$\rom; a $b$-linear mapping
$W\to V$ is nuclear if and only if it can be represented as
a sum \rom(possibly, infinite\rom)
of operators of rank~\rom1.
\endproclaim

The authors intend to carry out a detailed study of
$b$-nuclear operators in connection with idempotent analogs of
kernel theorems (in the spirit of Schwartz and Grothendieck)
in one of the subsequent papers. Here we only consider
the following example.

\example{Example~\rm(see~\cite{17, 18})}
Let $W = \sB(X,K)$ and $V = \sB(Y,K)$ (as described in~\S7).
By the kernel theorem~\cite{18},
every $b$-linear
operator $f\:W\to V$ has the form
$$
f\:\phi(x)\mapsto\wt\phi(y)=
 \int^\oplus_XK_f(x,y)\odot\phi(x)\,dx=
 \sup_{x\in X}(K_f(x,y)\odot\phi(x)),
\tag10
$$
where $\phi\in W$, $K_f(x,y)\in\sB(X\times Y,K)$ and the
``idempotent integral''
$\int^\oplus_X\psi(x)\,dx$
(e.g., see~\cite{2, 4--8})
is defined as $\sup_{x\in X}\psi(x)$ for each function
$\psi\in\sB(X,K)$. We define a $b$-linear
functional $K_{f,y}$ on $W$ by setting
$$
K_{f,y}(\phi) = \int^\oplus_XK_f(x,y)\phi(x)\,dx.
$$
The integral representation~\thetag{10} can be rewritten in
the form
$f(\phi) = \bigoplus_{y\in Y}K_{f,y}(\phi)\odot\delta_y$.
It follows from the kernel theorem that each $b$-linear
operator $f\:\sB(X,K)\to\sB(Y,K)$ is $b$-nuclear. It follows from
the theorem on the structure of a $b$-linear functional~\cite{18}
that each $b$-nuclear operator $f\:W\to V$
can be specified by an ``integral'' kernel $K_f(x,y)$.
Thus, the existence of an integral representation of an operator
is provided by the fact that the operator is nuclear and by
the existence of integral
representations of $b$-linear functionals. This scheme can be
generalized to a wide class
of $b$-complete semimodules.
\endexample

\remark{Remark~7}
Note that $W^* = (\wh W_b)^*$, where $\wh W_b$ is the
$b$-completion of the semimodule $W$ (the procedure of $b$-completion
for
semimodules is discussed, for example, in~\cite{12}).
Hence
the requirement that the semimodule $W$ be $b$-complete
is in fact not very important. In general, if semimodules $V$
and $W$
over $K$
admit $b$-completions $\wh V_b$ and $\wh W_b$ over $\wh K_b$,
then the ``completed'' $b$-tensor
product $V\wh\otimes W$ can be defined as
$\wh V_b\otimes\wh W_b$.
\endremark

\remark{Remark~8}
For the case of $a$-complete idempotent semimodules (and
tensors products of finite families of semimodules of this
type), some results of the present paper are also contained (in
other terms and under a different approach) in~\cite{19},
where a rather general categorical approach to tensor
products, suggested in~\cite{20}
(see also~\cite{21}), was studied.
Note that the approach of~\cite{20}
is not always convenient for applications to analysis in the
spirit of A. Grothendieck~\cite{1}.
For example, for the category of locally convex (or Banach)
spaces this approach provides only one of all possible (and
important in analysis)
topological tensor products in the sense of~\cite{1}.
In forthcoming papers, we shall consider various versions of
idempotent analogs of topological tensor
products.
\endremark

%





\centerline{\bf Acknowledgements}

This research was supported
by INTAS and the Russian Foundation for Basic Research
under  joint grant~No.~95-91.


\Refs


\item{1.}
A.~Grothendieck,
{\sl Produits tensoriels topologiques et espaces nucl\'eairs},
Mem. Amer. Math. Soc., {\bf 16},
Providence (R.I.)
(1955).

\item{2.}
V.~P.~Maslov,
{\sl Asymptotic Methods for Pseudodifferential Equations\/}
[in Russian],
Nauka,
Moscow
(1987).

\item{3.}
V.~P.~Maslov,
``A new superposition principle for optimization problems,''
{\sl Uspekhi Mat. Nauk\/},
{\bf 42},
No.~3,
39--48
(1987)
[in Russian; English Transl. in: Russian Math. Surveys {\bf 42}
(1987), no. 3, 43--54.

\item{4.}
V.~P.~Maslov,
{\sl M\'ethodes op\'eratorielles},
Mir,
Moscow
(1987).

\item{5.}
V.~P.~Maslov and S.~N.~Samborski\u\i, editors,
{\sl Idempotent Analysis},
Vol.~13 in {\sl Adv. Soviet Math},
Amer. Math. Soc.,
Providence (R.I.)
(1992).

\item{6.}
V.~P.~Maslov and V.~N.~Kolokoltsov,
{\sl Idempotent Analysis and Its Applications in Optimal
Control Theory\/}
[in Russian],
Nauka,
Moscow
(1994).

\item{7.}
V.~N.~Kolokoltsov and V.~P.~Maslov,
{\sl Idempotent Analysis and Applications},
Kluwer Acad. Publ.,
Dordrecht
(1997).

\item{8.}
G.~L.~Litvinov and V.~P.~Maslov,
{\sl Correspondence Principle for Idempotent Calculus and Some
Computer Applications},
Preprint IHES/M/95/33,
Institut des Hautes Etudes Scientifiques,
Bures-sur-Yvette
(1995) (see also~\cite{9, pp.~420--443}, and
{\tt http://arXiv.org/abs/math.GM/0101021}).

\item{9.}
{\sl Idempotency}
(J.~Gunawardena, editor),
Publ. of the Newton Institute,
Cambridge Univ. Press,
Cambridge
(1998).

\item{10.}
G.~Birkhoff,
{\sl Lattice Theory},
3d edition, American Mathematical Society Colloquium
Publications, Vol.~XXV,
American Mathematical Society, Providence, R.I. (1967).

\item{11.}
G.~L.~Litvinov, V.~P.~Maslov, and G.~B.~Shpiz,
``Linear functionals on idempotent spaces.
An algebraic approach,''
{\sl Dokl. Ross. Akad. Nauk\/},
{\bf 363}, No. 3, 298--300 (1998);
English transl. in: Doklady Mathematics, {\bf 58}, No. 3,
389--391 (1998) (see also {\tt http://arXiv.org/abs/math.FA/0012268}).

\item{12.}
G.~L.~Litvinov, V.~P.~Maslov, and G.~B.~Shpiz,
{\sl Idempotent Functional Analysis. I. An Algebraic
Approach\/},
Preprint,
International Sophus Lie Center,
Moscow
,1998 (see also {\tt http://arXiv.org/abs/math.FA/0009128}).

\item{13.}
F.~L.~Bacelli, G.~Cohen, G.~J.~Olsder, and J.~-P.~Quadrat,
{\sl Synchronization and Linearity: an Algebra for Discrete Event
Systems},
Wiley,
New York
(1992).

\item{14.}
M.~Gondran and M.~Minoux,
{\sl Graphes et algorithms},
Clarendon Paperbacks,
Oxford Univ. Press,
Oxford
(1979).

\item{15.}
H.~Sch\"afer,
{\sl Topological Vector Spaces},
   Graduate Texts in Mathematics, Vol.~3,
Springer-Verlag, New York--Berlin (1971).

\item{16.}
S.~I.~Gelfand and Yu.~I.~Manin,
{\sl Methods in Homological Algebra. Vol.~1. Introduction to
Cohomology Theory and Derived Categories\/}
[in Russian],
Nauka,
Moscow
(1988).

\item{17.}
P.~I.~Dudnikov and S.~N.~Samborskii,
``Endomorphisms of semimodules over semirings with idempotent
operations,''
(Kiev: Institute for Mathematics, Akad. Nauk Ukr. SSR, 1987),
{\sl Izv. Akad. Nauk SSSR Ser. Mat.} [{\sl Math. USSR-Izv.}],
{\bf 55},
No.~1,
91--105
(1991).

\item{18.}
M.~A.~Shubin,
{\sl Algebraic Remarks on Idempotent Semirings and a Kernel
Theorem in Spaces of Bounded Functions\/}
[in Russian],
Institute for New Technology,
Moscow
(1990).
English transl. in [5], 151--166.

\item{19.}
A.~Joyal and M.~Tierney,
``An extension of the Galois theory of Grothendieck,''
{\sl Mem. Amer. Math. Soc.},
{\bf 51},
No.~309
(1984).

\item{20.}
B.~Banaschewski and E.~Nelson,
``Tensor products and bimorphisms,''
{\sl Canad. Math. Bull.},
{\bf 19},
No.~4,
385--402
(1976).

\item{21.}
D.~Pumpl\"un,
``Das Tensorprodukt als universelles Problem,''
{\sl Math. Ann.},
{\bf 171},
247--262
(1967).

\endRefs




\bigskip

International Sophus Lie Center (G.~L.~Litvinov and G.~B.~Shpiz)

litvinov\@islc.msk.su

\bigskip

  Moscow State University (V.~P.~Maslov)

maslov\@ipmnet.ru


\enddocument